\documentclass[conference] {IEEEtran}
\usepackage[utf8]{inputenc}
\usepackage[reqno]{amsmath}
\usepackage{amssymb}
\usepackage{amsthm}
\usepackage{amscd}
\usepackage{enumerate}
\usepackage[table]{xcolor}
\usepackage{graphicx}

\theoremstyle{plain}

\newtheorem{lemma}{Lemma}[section]

\newtheorem{theorem}[lemma]{Theorem}

\theoremstyle{definition}

\theoremstyle{remark}

\numberwithin{equation}{section}

\begin{document}

\title{Exact Line Packings from Numerical Solutions}

\author{
\IEEEauthorblockN{Dustin G. {\sc Mixon} \& Hans {\sc Parshall}}

\IEEEauthorblockA{Department of Mathematics\\
The Ohio State University\\
Columbus, OH 43210\\
$\{${\tt mixon.23,parshall.6$\}$@osu.edu}}

%
}
\maketitle

\begin{abstract}
Recent progress in Zauner's conjecture has leveraged deep conjectures in algebraic number theory to promote numerical line packings to exact and verifiable solutions to the line packing problem.
We introduce a numerical-to-exact technique in the real setting that does not require such conjectures.
Our approach is completely reproducible, matching Sloane's database of putatively optimal numerical line packings with Mathematica's built-in implementation of cylindrical algebraic decomposition.
As a proof of concept, we promote a putatively optimal numerical packing of eight points in the real projective plane to an exact packing, whose optimality we establish in a forthcoming paper.
\end{abstract}

\section{Introduction}

Select $\mathbf{F}\in\{\mathbf{R},\mathbf{C}\}$, and consider the so-called line packing problem of packing $n$ points in $\mathbf{F}\mathbf{P}^{d-1}$ so that the minimum distance is maximized.
This fundamental problem resembles the classical Tammes problem~\cite{tammes30}, originally posed in 1930, which seeks to pack points on the sphere $S^2$, and indeed, the Tammes problem is equivalent to special case in which $\mathbf{F}=\mathbf{C}$ and $d=2$.
The general line packing problem was originally studied in the 1960s and 70s by Fejes T\'{o}th~\cite{fejes65}, Welch~\cite{welch74}, Delsarte, Goethals and Seidel~\cite{delsarte75}, and Levenshtein~\cite{levenshtein82}.
In 1996, Conway, Hardin and Sloane~\cite{conway96} rejuvenated interest in this problem, providing a plethora of putatively optimal packings (available in Sloane's online database~\cite{sloane}) and also proving the so-called orthoplex bound.
The last decade of research has identified various applications of optimal line packings, including compressed sensing~\cite{bandeira17}, digital fingerprinting~\cite{mixon13}, quantum state tomography~\cite{renes04}, and multiple description coding~\cite{strohmer03}.
This in turn has sparked a flurry of work to construct optimal packings.
Most of this work finds new packings that achieve equality in the Welch bound (see~\cite{fickus16} for a survey), though there has also been progress in achieving equality in the orthoplex bound~\cite{bodmann16} and the Levenshtein bound~\cite{haas17}.
In addition, last year, Bukh and Cox~\cite{bukh18} discovered a new bound along with corresponding optimal packings.
Despite the substantial progress, optimal packings remain unidentified for the vast majority of triples $(\mathbf{F},d,n)$.
In terms of real degrees of freedom, the smallest open case to date is $(\mathbf{R},3,8)$.

The application of quantum state tomography concerns a rather interesting instance of the complex case.
In particular, Zauner~\cite{zauner99} conjectured that for every $d>1$, the optimal packings of $d^2$ points in $\mathbf{CP}^{d-1}$ necessarily achieve equality in the Welch bound, and furthermore, they can be constructed to exhibit symmetries from the Weyl--Heisenberg group.
After a decade of significant effort to prove Zauner's conjecture~\cite{fuchs17}, the conjecture is only known to hold for finitely many choices of $d$.
Some of the latest work along these lines has leveraged an observation that all of the known optimal packings of Zauner's form feature coordinates that reside in a predictable number field~\cite{appleby18}.
One may exploit this observation to promote a numerical solution to an exact solution: obtain thousands of digits of precision, use these digits to guess the exact coordinates, and then use symbolic calculations to verify that the resulting packing achieves equality in the Welch bound.
This procedure has led to several new constructions.

In the present paper, we take inspiration from~\cite{appleby18} to develop a completely different technique that promotes numerical packings to exact packings.
This approach is made possible by Sloane's database~\cite{sloane}, which provides many putatively optimal numerical packings in the real case.
Unlike the Zauner instance, we do not have access to conjectures that predict the field structure of optimal packing coordinates.
Instead, we will borrow ideas introduced in~\cite{fickus18}, which leveraged quantifier elimination over the reals to find computer-assisted proofs of certain optimal line packings.
We will focus on the case of 8 points in $\mathbf{RP}^2$ (i.e., the smallest open case to date), but our methods generalize to arbitrary real packings.
Sadly, this exact packing does not achieve equality in any known bound, and so we cannot simply verify optimality by symbolic calculation, as in the Zauner instance.
Instead, we prove optimality in a forthcoming paper using a computational graph theory approach that takes inspiration from recent progress on the Tammes problem~\cite{musin12,musin15}.

In the next section, we set notation and review the necessary background before stating the main result.
Section~\ref{proofsection} then outlines the computer-assisted proof of our result, which we document carefully for the sake of reproducibility.
We conclude in Section~\ref{futuresection} with various directions for future work.

\section{Main Result}

\begin{figure*}
\[
G_\mathsf{CHS}:=
\left(\begin{array}{rrrrrrrr}
1.00000 & -0.59840 & \cellcolor{lightgray}{0.64759} & -0.12425 & 0.16026 & \cellcolor{lightgray}{-0.64759} & -0.22283 & \cellcolor{lightgray}{-0.64759} \\
-0.59840 & 1.00000 & \cellcolor{lightgray}{-0.64759} & -0.44579 & \cellcolor{lightgray}{-0.64759} & -0.01815 & \cellcolor{lightgray}{-0.64759} & 0.39325 \\
\cellcolor{lightgray}{0.64759} & \cellcolor{lightgray}{-0.64759} & 1.00000 & \cellcolor{lightgray}{0.64759} & -0.16026 & \cellcolor{lightgray}{-0.64759} & 0.16026 & 0.10359 \\
-0.12425 & -0.44579 & \cellcolor{lightgray}{0.64759} & 1.00000 & -0.10359 & -0.01815 & \cellcolor{lightgray}{0.64759} & \cellcolor{lightgray}{0.64759} \\
0.16026 & \cellcolor{lightgray}{-0.64759} & -0.16026 & -0.10359 & 1.00000 & \cellcolor{lightgray}{0.64759} & \cellcolor{lightgray}{0.64759} & \cellcolor{lightgray}{-0.64759} \\
\cellcolor{lightgray}{-0.64759} & -0.01815 & \cellcolor{lightgray}{-0.64759} & -0.01815 & \cellcolor{lightgray}{0.64759} & 1.00000 & \cellcolor{lightgray}{0.64759} & -0.01815 \\
-0.22283 & \cellcolor{lightgray}{-0.64759} & 0.16026 & \cellcolor{lightgray}{0.64759} & \cellcolor{lightgray}{0.64759} & \cellcolor{lightgray}{0.64759} & 1.00000 & 0.12425 \\
\cellcolor{lightgray}{-0.64759} & 0.39325 & 0.10359 & \cellcolor{lightgray}{0.64759} & \cellcolor{lightgray}{-0.64759} & -0.01815 & 0.12425 & 1.00000
\end{array}\right)
\]
\caption{The Gram matrix of the putatively optimal 8-packing in $\mathbf{RP}^2$, from Sloane's online database~\cite{sloane}.}
\label{sloaneGram}
\end{figure*}

We define an $n$-packing in $\mathbf{RP}^{d - 1}$ to be a sequence of $n$ lines through the origin in $\mathbf{R}^d$, which we identify with a $d \times n$ matrix $\Phi = [\varphi_1 \cdots  \varphi_n]$ whose columns are unit vectors spanning the corresponding lines.  Define the coherence by
\[
	\mu(\Phi) := \max_{1 \leq i < j \leq n} | \langle \varphi_i, \varphi_j \rangle |.
\]
We say $\Phi$ is optimal when $\mu(\Phi) \leq \mu(\Phi')$ for every $n$-packing $\Phi'$ in $\mathbf{RP}^{d - 1}$.

When searching for an optimal $n$-packing, it is equivalent to search for a corresponding Gram matrix $\Phi^T \Phi$.  A real symmetric $n \times n$ matrix $G$ is the Gram matrix for some $n$-packing $\Phi$ in $\mathbf{RP}^{d - 1}$ if and only if $\operatorname{rank}(G) \leq d$, $G$ is positive semidefinite, and its diagonal entries each satisfy $G_{ii} = 1$.  From this perspective, the largest off-diagonal entries of $G$ achieve the coherence:
\[
	\max_{1 \leq i < j \leq n} |G_{ij}| = \mu(\Phi).
\]

The explicit numerical solutions of Conway, Hardin and Sloane~\cite{conway96} provide upper bounds on the coherence of optimal $n$-packings.  For the case of an $8$-packing in $\mathbf{RP}^2$, their corresponding Gram matrix $G_{\mathsf{CHS}}$ is depicted in Fig.~\ref{sloaneGram}, where we have rounded the corresponding entries for display on the page.  The Conway--Hardin--Sloane packing witnesses that any optimal $8$-packing $\Phi$ in $\mathbf{RP}^2$ necessarily satisfies $\mu(\Phi) \leq 0.64759$.

Several lower bounds on coherence are known, but the bound most relevant to $8$-packings in $\mathbf{RP}^2$ is originally due to Levenshtein~\cite{levenshtein82}; see~\cite{haas17} for a recent account.

\begin{lemma}
	If $\Phi$ is an $n$-packing in $\mathbf{RP}^{d - 1}$, then
	\[
		\mu(\Phi) \geq \sqrt{\frac{3n - d^2 - 2d}{(n - d)(d + 2)}}.
	\]
	In particular, if $\Phi$ is an $8$-packing in $\mathbf{RP}^2$, then $\mu(\Phi) \geq 0.6$.
\end{lemma}

While many of the entries in the matrix in Fig.~\ref{sloaneGram} appear to be equal up to numerical precision, we have no guarantee that these relationships must hold.  Regardless, we proceed by assuming that the entries that are roughly $\pm 0.64759$ are indeed equal in absolute value. To justify this assumption, we consider the contact graph of an $n$-packing $\Phi = [\varphi_1 \cdots \varphi_n]$, where the vertices are the numbers $\{1, \ldots, n\}$ and we connect $i$ to $j$ with an edge exactly when $|\langle \varphi_i, \varphi_j \rangle| = \mu(\Phi)$.  In forthcoming work, we prove that the adjacency matrix of the contact graph of the optimal 8-packing in $\mathbf{RP}^2$ can be obtained from $G_{\mathsf{CHS}}$ by replacing the entries of $\pm 0.64759$ with 1 and zeroing out all other entries.  See Fig.~\ref{graphFig} for the resulting embedding of this graph in the projective plane.

What follows is our main result:

\begin{figure}
\centering
\setlength{\unitlength}{0.35\textwidth}
\begin{picture}(1,1)
\put(0,0){\includegraphics[width=0.35\textwidth]{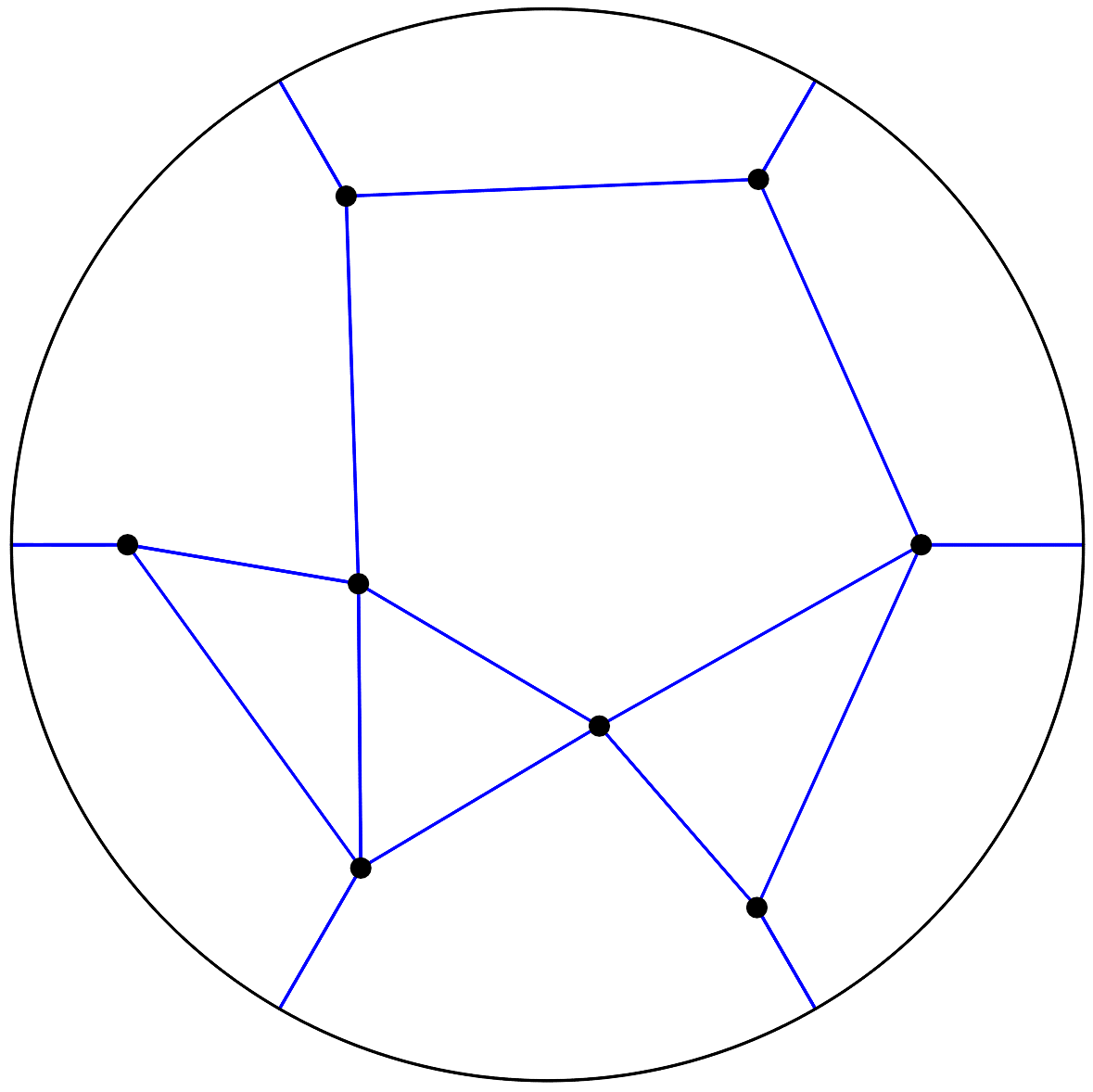}}
\put(0.64,0.15){$1$}
\put(0.105,0.525){$2$}
\put(0.785,0.5){$3$}
\put(0.67,0.775){$4$}
\put(0.35,0.47){$5$}
\put(0.54,0.36){$6$}
\put(0.35,0.16){$7$}
\put(0.335,0.76){$8$}
\end{picture}
\caption{The projective planar embedding of the contact graph of the packing corresponding to $G_{\mathsf{CHS}}$.  For instance, vertices $1$ and $8$ are connected by an edge.  Vertex labels correspond to the order of the packing as given in~\cite{sloane}.}
\label{graphFig}
\end{figure}

\begin{theorem}\label{main}

There is a unique real matrix of the form
\[
G = \left(\begin{array}{rrrrrrrr}
1 & a_1 & \cellcolor{lightgray}{\mu} & a_2 & a_3 & \cellcolor{lightgray}{-\mu} & a_4 & \cellcolor{lightgray}{-\mu} \\
a_1 & 1 & \cellcolor{lightgray}{-\mu} & a_5 & \cellcolor{lightgray}{-\mu} & a_6 & \cellcolor{lightgray}{-\mu} & a_7 \\
\cellcolor{lightgray}{\mu} & \cellcolor{lightgray}{-\mu} & 1 & \cellcolor{lightgray}{\mu} & a_8 & \cellcolor{lightgray}{-\mu} & a_9 & a_{10} \\
a_2 & a_5 & \cellcolor{lightgray}{\mu} & 1 & a_{11} & a_{12} & \cellcolor{lightgray}{\mu} & \cellcolor{lightgray}{\mu} \\
a_3 & \cellcolor{lightgray}{-\mu} & a_8 & a_{11} & 1 & \cellcolor{lightgray}{\mu} & \cellcolor{lightgray}{\mu} & \cellcolor{lightgray}{-\mu} \\
\cellcolor{lightgray}{-\mu} & a_6 & \cellcolor{lightgray}{-\mu} & a_{12} & \cellcolor{lightgray}{\mu} & 1 & \cellcolor{lightgray}{\mu} & a_{13} \\
a_4 & \cellcolor{lightgray}{-\mu} & a_9 & \cellcolor{lightgray}{\mu} & \cellcolor{lightgray}{\mu} & \cellcolor{lightgray}{\mu} & 1 & a_{14} \\
\cellcolor{lightgray}{-\mu} & a_7 & a_{10} & \cellcolor{lightgray}{\mu} & \cellcolor{lightgray}{-\mu} & a_{13} & a_{14} & 1
\end{array}\right)
\]
with the following properties:
\begin{itemize}
	\item[(i)] $\operatorname{rank}(G) = 3$,
	\item[(ii)] $G$ is positive semidefinite,
	\item[(iii)] $0.6 \leq \mu \leq 0.64759$, and
	\item[(iv)] $|a_j| < \mu$ for all $1 \leq j \leq 14$.
\end{itemize}
Moreover $\mu = \mu_0$, where $\mu_0$ is the largest root of
\[
	\scriptstyle 1 + 5x - 8x^2 - 80x^3 - 78x^4 + 146x^5 - 80x^6 - 584x^7 + 677x^8 + 1537x^9
\]
and is given numerically by $\mu_0 \approx 0.6475889787$.
\end{theorem}

The matrix in Theorem~\ref{main} equals $G_{\mathsf{CHS}}$ up to precision.

\begin{table*}
\centering

\caption{Outline of CAD-assisted proof of $\mu = \mu_0$}
\begin{tabular}{c | l | l | l | l}
\hline
Step & Rows & Columns & Variable Order & Constraints Obtained\\
\hline
1 & 2, 5, 6, 7 & 2, 5, 6, 7 & $\mu, a_6$ & $a_6 = (1 + \mu - 4\mu^2) / (1 + \mu)$\\

2 & 2, 3, 5, 6, 7 & 2, 3, 5, 6, 7 & $\mu, a_8, a_9$ & $a_8^2 = (1 + \mu - 3\mu^2 - \mu^3) / (2 + 4\mu)$; $a_9 = -a_8$\\

3 & 1, 2, 3, 5, 6, 7 & 1, 2, 3, 5, 6, 7 & $\mu, a_8, a_3, a_1, a_4$ & $a_4 = (-a_1 - \mu - a_1\mu - 2a_3\mu - \mu^2) / (2\mu)$; $\mu \geq (1 + \sqrt{17})/8 \approx 0.64038$\\

4 & 3, 4, 6, 7 & 3, 4, 6, 7 & $\mu, a_8, a_{12}$ & $a_{12}^2 = a_6^2$\\

5 & 2, 3, 4, 6, 7 & 2, 4, 6, 7 & $\mu, a_5, a_{12}, a_8$ & $a_{12} = a_6$; $a_8 < 0$\\

6 & 2, 3, 5, 6 & 1, 2, 3, 5, 6, 7 & $\mu, a_8, a_3, a_1$ & $a_1 =$ $\scriptstyle  (-\mu - a_8 \mu - 4 a_3 a_8 \mu + 2 a_8^2 \mu + 3 \mu^2 - a_8 \mu^2 - 
 2 \mu^3) / (-1 + a_8 + 2 a_8^2 + \mu + a_8 \mu)$\\

7 & 1, 2, 3, 5, 6 & 1, 2, 3, 5, 6 & $\mu, a_3, a_8$ & $a_3 = -a_8$\\

8 & 4, 5, 6, 7 & 2, 3, 5, 6, 7 & $\mu, a_8, a_5, a_{11}$ & $a_{11} = (1 - a_5 + \mu - a_5 \mu - 6 \mu^2) / (2 \mu)$\\

9 & 1, 3, 4, 6, 7 & 1, 3, 4, 6, 7 & $\mu$, $a_2$, $a_8$ & $a_2 = (-3\mu - 2\mu^2 + 9\mu^3) / (1 + 2\mu + \mu^2)$\\

10 & 1, 4, 5, 6, 8 & 1, 4, 5, 6, 8 & $\mu$, $a_8$, $a_5$, $a_{13}$ & $\mu = \mu_0$\\
\hline
\end{tabular}

\label{prooftable}

\end{table*}

\section{Proof of Main Result} \label{proofsection}

In principle, Theorem~\ref{main} amounts to quantifier elimination over the reals.  Indeed, we can ensure that $\operatorname{rank}(G) \leq 3$ by ensuring that each of its $4 \times 4$ minors vanish. Furthermore, by Sylvester's criterion, $G$ is positive semidefinite if and only if all of its principal minors are nonnegative.  Hence, the admissible choices for $G$ are precisely the solutions to a finite collection of polynomial equalities and inequalities.  For relatively small problems, this can be accomplished somewhat efficiently by appealing to the cylindrical algebraic decomposition (CAD) algorithm introduced by Collins~\cite{collins75}, for which we use the implementation available in Mathematica.  The main idea behind CAD involves constructing a projection from our solution set to a semialgebraic set of one dimension lower, eliminating a variable.  This process is iterated until a subset of $\mathbf{R}$ is reached, at which point the desired semialgebraic decomposition of the solution set can be obtained by iteratively lifting.  For a useful introduction to using CAD with Mathematica, see~\cite{kauers10}.

Unfortunately, the runtime of CAD is doubly exponential in the number of variables, which was shown to be intrinsic to the problem by Davenport and Heintz~\cite{davenport88}.  Moreover, the speed of the algorithm is highly sensitive to the order in which the projections are constructed, i.e, the order in which the variables are specified.  It may come as no surprise that for our 15-variable system, naive CAD queries fail to deliver the desired decomposition in a reasonable amount of time.  As such, we proceed by iteratively applying CAD to subsystems of the original problem in order to successively reduce the number of variables.

In Table~\ref{prooftable}, we outline our computer-assisted proof that $\mu = \mu_0$.  In each step, we select a number of rows and columns to analyze, insisting that each $4 \times 4$ minor vanishes and the bounds on $\mu$ and each $|a_j|$ are simultaneously satisfied.  In practice, we found that it was more efficient to ignore the positive semidefinite constraint and verify at the end that our solution for $G$ satisfied Sylvester's criterion.  For each CAD query, we must also choose an ordering of the variables involved, and CAD reports the first listed variable as independent and each subsequent variable depending upon possibly all of the preceding variables.    Our motivation for each selection was to keep the number of variables small so that CAD would terminate quickly, and overall, our approach computes $\mu = \mu_0$ in ten steps that take a total of roughly 5 minutes in Mathematica 11 on a 3.4 GHz Intel Core i5.

In Step~1 of our procedure, we use the rows and columns with indices in $\{2, 5, 6, 7\}$, since the resulting minor produces the two-variable equation
\[
\left|
\begin{array}{rrrr}
1 & -\mu & a_6 & -\mu \\ -\mu & 1 & \phantom{-}\mu & \mu \\ a_6 & \mu & 1 & \mu \\ -\mu & \mu & \mu & 1
\end{array}\right| = 0,
\]
which amounts to
\[
(a_6 + 1)(\mu - 1)(a_6(\mu + 1) + 4\mu^2 - \mu - 1) = 0.
\]
The restrictions $|a_6| < \mu \leq 0.64759$ then guarantee that
\[
a_6 = \frac{1 + \mu - 4\mu^2}{1 + \mu}.
\]
This allows us to replace $a_6$ with a rational function of $\mu$ in all subsequent computations.

Next, in Step~2, we consider the $4 \times 4$ minors with rows and columns $\{2, 3, 5, 6,7\}$ to obtain five polynomial equations in the three variables $\mu, a_8, a_9$.  Notice that we would have been burdened by an additional variable had we not already performed Step~1 to eliminate $a_6$.  These equations lead CAD to report $a_9 = -a_8$, as one might have guessed from inspecting Fig.~\ref{sloaneGram}. Importantly, this reduces the number of variables once again. CAD also reports the new constraint
\[
	a_8 = \pm \sqrt{\frac{1 + \mu - 3\mu^2 - \mu^3}{2 + 4\mu}}.
\]
Even without the sign ambiguity, which is eventually removed in Step~5, this does not quite allow us to replace the variable $a_8$ with a rational function in $\mu$.  For this reason, we continue treating $a_8$ as a free variable for the remainder of our proof.

Steps~3--10 proceed in a similar manner, collecting constraints until finally arriving at $\mu = \mu_0$.  At this point, CAD reports exact expressions for every coefficient of $G$ except $a_7$, $a_{10}$, and $a_{14}$.  These remaining coefficients can be determined from the following principal minors: $\{2,5,6,8\}$ determines~$a_7$, $\{3,4,6,8\}$ determines $a_{10}$, and $\{5,6,7,8\}$ determines $a_{14}$.  Each coefficient of $G$ has algebraic degree 9, and the entries apparently equal in absolute value in $G_{\mathsf{GHS}}$ are indeed so.

We still need to ensure that the resulting matrix $G$ has rank 3 and is positive semidefinite.  Mathematica quickly reports $\operatorname{rank}(G) = 3$, since this only requires verifying that a single $3 \times 3$ minor is non-vanishing.  However, querying whether $G$ is positive semidefinite does not terminate in a reasonable amount of time, since computing either its eigenvalues or its $3 \times 3$ principal minors symbolically is computationally expensive.  To avoid this, we simply apply a perturbation argument: it is more than enough to ask Mathematica to report a numerical approximation of $G$ to within 100 digits of precision and then numerically verify that each $3 \times 3$ principal minor of the approximation is bounded below by $0.0001$.

\section{Future Work} \label{futuresection}

In this paper, we provide a new approach for finding exact line packings from numerical solutions. We focused on the special case of $8$ points in $\mathbf{RP}^2$, but there is hope to apply our method to many more instances in the real case.  Sloane~\cite{sloane} provides putatively optimal numerical $n$-packings for $n \leq 100$ in $\mathbf{RP}^2$ (as well as various numerical packings in $\mathbf{RP}^{d-1}$ for $d\leq 16$), and we suspect that CAD can be applied productively well beyond our test case.  The next case of $9$ points in $\mathbf{RP}^2$ is believed to have two essentially distinct optimal line packings~\cite{conway96}.  Solving for these packings exactly would allow us to verify that their coherences are indeed equal.

It would be particularly useful to automate our seemingly ad-hoc choices of efficiently solvable subsystems in Table~\ref{prooftable}.  For example, the so-called sketch-and-solve paradigm of solving random subsystems has found success in solving large least-squares systems~\cite{woodruff14}.  Recently, Huang~et.~al~\cite{huang14} applied machine learning to select between various popular heuristics for ordering the variables in a CAD query.  Presumably, one can effectively apply machine learning to assist in the construction of exact line packings.

Alternatives to CAD could plausibly speed up our proof in Table~\ref{prooftable}, and more interestingly, allow for the computation of even more exact packings.  Indeed, this numerical-to-exact approach could serve as a benchmark for various alternatives to CAD that may emerge from the computational algebraic geometry community.

Since our method relies on quantifier elimination over the reals, applying it to the complex case necessarily involves doubling the number of variables.  This presents another reason to pursue alternatives to CAD.  Also, Sloane's database only provides numerical packings in the real case, and so we desire a similar table for the complex case.

While we have computed an exact version of the packing corresponding to $G_{\mathsf{CHS}}$, it remains to demonstrate that this is in fact an optimal packing, as conjectured in~\cite{conway96}.  Most proofs of optimality in the literature proceed by establishing equality in a general lower bound for coherence.  This approach seems to have limitations, since most of the known numerical solutions for the line packing problem do not approach this threshold.  In our forthcoming proof of optimality, we use computational graph theory methods that do not yet have an appropriate analogue in $\mathbf{RP}^{d - 1}$ for $d > 2$.  In particular, we make use of projective planarity, much like how recent computational solutions to the Tammes problem leverage planarity~\cite{musin12,musin15}.

\section*{Acknowledgments}

DGM was partially supported by AFOSR FA9550-18-1-0107, NSF DMS 1829955, and the Simons Institute of the Theory of Computing.


\begin{thebibliography}{11}

\bibitem{appleby18}
\textsc{M. Appleby, T. Y. Chien, S. Flammia \& S. Waldron}
\newblock{\em Constructing exact symmetric informationally complete measurements from
numerical solutions.}
{J. Phys. A} {\bf 51} (2018), 165302.

\bibitem{bandeira17}
\textsc{A. S. Bandeira, M. Fickus, D. G. Mixon \& P. Wong}
\newblock{\em The road to deterministic matrices with
the restricted isometry property.}
{J. Fourier Anal. Appl.} {\bf 19} (2013), 1123-1149.

\bibitem{bodmann16}
\textsc{B. G. Bodmann \& J. Haas}
\newblock{\em Achieving the orthoplex bound and constructing weighted complex projective 2-designs with Singer sets.}
{Linear Algebra Appl.} {\bf 511} (2016) 54--71.

\bibitem{bukh18}
\textsc{B. Bukh \& C. Cox}
\newblock{\em Nearly orthogonal vectors and small antipodal spherical codes.}
preprint, arXiv:1803.02949, (2018).

\bibitem{collins75}
\textsc{G. E. Collins}
\newblock{\em Quantifier elimination for real closed fields by cylindrical algebraic decomposition.}
{Automata theory and formal languages, Second GI Conf., Kaiserslautern} (1975), 134--183.

\bibitem{conway96}
\textsc{J. H. Conway, R.H. Hardin \& N. J. A. Sloane}
\newblock{\em Packing Lines, Planes, etc.: Packings in Grassmannian Spaces.}
{Exp. Math.} {\bf 5} (1996), 139--159.

\bibitem{davenport88}
\textsc{J. H. Davenport \& J. Heintz}
\newblock{\em Real quantifier elimination is doubly exponential.}
{Journal of Symbolic Computation} {\bf 5} (1988), 29--35.

\bibitem{delsarte75}
\textsc{P. Delsarte, J. M. Goethals, J. J. Seidel}
\newblock{\em Bounds for systems of lines, and Jacobi polynomials.}
{Philips Res. Rep.} {\bf 30} (1975)

\bibitem{fejes65}
\textsc{L. Fejes T\'{o}th}
\newblock{\em Distribution of points in the elliptic plane.}
{Acta Math. Acad. Sci. Hungar} {\bf 16} (1965), 437--440.

\bibitem{fickus16}
\textsc{M. Fickus \& D. G. Mixon}
\newblock{\em Tables of the Existence of Equiangular Tight Frames.}
{preprint}, arXiv:
1504.00253, (2016).

\bibitem{fickus18}
\textsc{M. Fickus, J. Jasper \& D. G. Mixon}
\newblock{\em Packings in real projective spaces.}
{SIAM J. Appl. Algebra Geom.} {\bf 2} (2018) 377--409.

\bibitem{fuchs17}
\textsc{C. A. Fuchs, M. C. Hoang \& B. C. Stacey}
\newblock{\em The SIC Question: History and State of Play.}
{Axioms.} {\bf 6} (2017) 21.

\bibitem{haas17}
\textsc{J. I. Haas, N. Hammen \& D. G. Mixon}
\newblock{\em The Levenstein bound for packings in projective spaces.}
{Wavelets and Sparsity XVII} (2017), 103940V.

\bibitem{huang14}
\textsc{Z. Huang, M. England, D. Wilson, J. H. Davenport, L. Paulson \& J. Bridge}
\newblock{\em Applying machine learning to the problem of choosing a heuristic to select the variable
ordering for cylindrical algebraic decomposition.}
{Intelligent Computer Mathematics} (2014) 92--107.

\bibitem{kauers10}
\textsc{M. Kauers}
\newblock{\em How to use cylindrical algebraic decomposition.}
{S\'{e}m. Lothar. Combin.} {\bf 65}  (2010/12), B65a.

\bibitem{levenshtein82}
\textsc{V. I. Levenshtein}
\newblock{\em Bounds on the maximal cardinality of a code with bounded modulus of the inner product.}
{Soviet Math. Dokl.} {\bf 25} (1982), 526--531.

\bibitem{mixon13}
\textsc{D. G. Mixon, C. J. Quinn, N. Niyavash \& M. Fickus}
\newblock{\em Fingerprinting with equiangular tight frames.}
{IEEE Trans. Inform. Theory} {\bf 59} (2013), 1855--1865.

\bibitem{musin12}
\textsc{O. R. Musin \& A. S. Tarasov}
\newblock{\em The Strong Thirteen Spheres Problem.}
{Discr. Comput. Geom.} {\bf 48} (2012), 128--141.

\bibitem{musin15}
\textsc{O. R. Musin \& A. S. Tarasov}
\newblock{\em The Tammes problem for $N = 14$.}
{Exp. Math.} {\bf 24} (2015), 460--468.

\bibitem{renes04}
\textsc{J. M. Renes, R. Blume-Kohout, A. J. Scott \& C. M. Caves}
\newblock{\em Symmetric informationally complete
quantum measurements} J. Math. Phys. {\bf 45} (2004), 2171--2180.

\bibitem{sloane}
\textsc{N. J. A. Sloane}
\newblock{\em Packings in Grassmannian spaces.}
{http://neilsloane.com/grass/}

\bibitem{strohmer03}
\textsc{T. Strohmer \& R. W. Heath}
\newblock{\em Grassmannian frames with applications to coding and communication.}
{Appl. Comput. Harmon. Anal.} {\bf 14} (2003), 257--275.

\bibitem{tammes30}
\textsc{R. M. L. Tammes}
\newblock{\em On the Origin Number and Arrangement of the Places of Exits on the Surface of Pollengrains.}
{Rec. Trv. Bot. Neerl.} {\bf 27} (1930) 1--84.

\bibitem{welch74}
\textsc{L. R. Welch}
\newblock{\em Lower bounds on the maximum cross correlation of signals.}
{IEEE Trans. Inform. Theory} {\bf 20} (1974), 397--399.

\bibitem{woodruff14}
\textsc{D. P. Woodruff}
\newblock{\em Sketching as a tool for numerical linear algebra.}
{Found. Trends Theor. Comput. Sci.} {\bf 10} (2014), 1--157.

\bibitem{zauner99}
\textsc{G. Zauner}
\newblock{\em Quantendesigns - Grundz\"{u}ge einer nichtkommutativen Designtheorie.}
{Ph.D. thesis, University of Vienna, Vienna, Austria, 1999.}


\end{thebibliography}
\end{document}